\documentclass[11pt]{combine}
\usepackage{amsmath,amsfonts,amssymb,amsthm,epsfig, graphicx, hyperref}
\usepackage[dvipsnames,table,dvipsnames*, svgnames*, hyperref]{xcolor}
\usepackage{natbib,amsmath,amssymb,amsthm,graphicx,setspace,paralist,booktabs,rotating,subcaption,float,color}
\usepackage{multirow}
\usepackage{fancyhdr}
\usepackage{bibunits,cases}
\usepackage[small]{titlesec}

\newtheorem{theorem}{Theorem}
\newtheorem{corollary}{Corollary}
\newtheorem{proposition}{Proposition}

\newtheorem{lemma}{Lemma}
{
\theoremstyle{definition}
\newtheorem{definition}{Definition}
\newtheorem{example}{Example}
\newtheorem{remark}{Remark}
}
\newcommand{\beq}{\begin{equation}}
\newcommand{\eeq}{\end{equation}}
\newcommand{\beas}{\begin{align*}}
\newcommand{\eeas}{\end{align*}}
\newcommand{\bea}{\begin{align}}
\newcommand{\eea}{\end{align}}
\newcommand{\bei}{\begin{itemize}}
\newcommand{\eei}{\end{itemize}}
\newcommand{\ben}{\begin{enumerate}}
\newcommand{\een}{\end{enumerate}}
\newcommand{\bet}{\begin{theorem}}
\newcommand{\eet}{\end{theorem}}
\newcommand{\bel}{\begin{lemma}}
\newcommand{\eel}{\end{lemma}}
\newcommand{\bep}{\begin{proposition}}
\newcommand{\eep}{\end{proposition}}

\newcommand{\bed}{\begin{definition}}
\newcommand{\eed}{\end{definition}}
\newcommand{\bec}{\begin{corollary}}
\newcommand{\eec}{\end{corollary}}
\newcommand{\bex}{\begin{example}}
\newcommand{\eex}{\end{example}}

\newcommand{\R}{\mathbb{R}}

\theoremstyle{plain}
\newtheorem{thm}{Theorem}[section]

\newtheorem{lem1}{Lemma A}

\newcommand{\f}{\frac}

\newcommand{\nn}{\nonumber}
\newcommand{\expit}{\text{expit}}
\newcommand{\logit}{\text{logit}}
\numberwithin{equation}{section}

\newcommand{\T}[1]{{\text{T}}}

\usepackage{algorithm,algorithmic}

\begin{document}


\title{\scshape On the global identifiability of logistic regression models with misclassified outcomes}
\author{Rui Duan$^1$, Yang Ning$^2$, Jiasheng Shi$^3$, Raymond J Carroll$^4$, Tianxi Cai$^1$ and Yong Chen$^5$ \\
$^1$ Department of Biostatistics, Harvard T.H. Chan School of Public Health, \\ Boston, Massachusetts 02115, U.S.A.\\
$^2$ Department of Statistics and Data Science, Cornell University,\\ Ithaca, New York 14853, U.S.A.\\
$^3$ Department of Research, The Children’s Hospital of Philadelphia,\\ Philadelphia, Pennsylvania 19104, U.S.A.\\
$^4$ Department of Statistics,  Texas A\&M University,\\ College Station, Texas 77843, U.S.A.\\
$^5$Department of Biostatistics, Epidemiology and Informatics, \\
University of Pennsylvania,
Philadelphia, PA 19104\\
}
\date{}
\maketitle
\thispagestyle{empty}

	\begin{abstract}
In the last decade, the secondary use of large data from health systems, such as electronic health records, has demonstrated great promise in advancing biomedical discoveries and improving clinical decision making. However, there is an increasing concern about biases in association studies caused by misclassification in the binary outcomes derived from electronic health records. We revisit the classical logistic  regression model with misclassified outcomes. 
Despite that local identification conditions in some related settings have been previously established, the global identification of such models remains largely unknown and is an important question yet to be answered. We derive necessary and sufficient conditions for global identifiability of logistic regression models with misclassified outcomes, using a novel approach termed as the submodel analysis, and a technique adapted from the Picard-Lindel\"{o}f existence theorem in ordinary differential equations. In particular, our results are applicable to logistic models with discrete covariates, which is a common situation in biomedical studies, 
The conditions are easy to verify in practice. In addition to model identifiability, we propose a hypothesis testing procedure for regression coefficients in the misclassified logistic regression model when the model is not identifiable under the null.

	\bigskip
	
\end{abstract}
 \linespread{1.8}
\section{Introduction}

The logistic regression model is one of the most widely used tools for studying associations between a binary outcome and multivariate covariates. It is ubiquitous in a variety of research areas including social science and  biomedical studies, due to the interpretation of log odds ratio for the regression coefficients \citep{hosmer2004applied} and various advantages in statistical inference, such as the availability of conditional inference \citep{lindsay1982conditional,stefanski1987conditional} and applicability in retrospective studies \cite{prentice1979logistic}. 

In some applications, the binary outcome may be subject to misclassification. For example, disease status ascertained by an imperfect diagnostic test. As a consequence, the estimated associations based on the error-prone outcomes, hereafter referred to as surrogate outcomes, are often biased, sometimes badly so. One example rises from the secondary use of electronic health records data, where a major concern of the use of electronic health records data is the bias caused by misclassification in the binary outcomes (phenotypes) derived from the data. Although manual chart review can be conducted to obtain the true phenotype of patients in some situations, such a process is oftentimes consuming and expensive. As a consequence, methodologies for reducing the number of validation samples or conducting the association analysis without the validation samples have been widely considered in recent work \citep{yu2017enabling, gronsbell2019automated}.

There is a literature studying the impacts of misclassification on the biases and statistical methods to correct for such biases. With internal or external validation samples, the impacts of misclassification can be characterized by comparing the estimated associations based on true outcomes with the estimated associations based on surrogate outcomes \citep{neuhaus99,duan2016empirical}, and the biases in estimated associations can be corrected through various types of calibration \citep{carroll2006measurement}. In other situations, validation samples can be difficult to obtain due to constraints on cost or time, ethical considerations, or invasiveness of the validation procedure. Under such situations, bias correction can be conducted via postulating a parametric assumption on the relationship between the underlying binary outcome and covariates, and a parametric characterization of the misclassification rates of the surrogate outcomes \citep{neuhaus99,neuhaus2002analysis,lyles2002note,lyles2010sensitivity}. However, as we elaborate below, model identifiability, which is fundamental for the applicability of standard statistical inference procedures, has not been fully investigated. There are at least two distinct definitions of the model identifiability: global and local identifiability. Suppose a random variable $Y$ has a probability density (or mass) function $f(y;\theta)$ indexed by parameter $\theta$. The parameter $\theta$ is said to be globally identifiable if $f(y; \theta^*)=f(y; \theta)$ for all $y$ almost surely implies $\theta^*=\theta$. This definition excludes the undesirable situation where different parameter values give rise to the same distribution, which makes it impossible to distinguish between these parameter values based on the data. A weaker form of identifiability is local identifiability. We say a model is locally identifiable at $\theta_0$ if there exists a neighborhood $\mathcal{B}$ of $\theta_0$ such that no other $\theta\in \mathcal{B}$ has the same probability distribution $f(y; \theta_0)$ as for all $y$ almost surely \citep{huang2004building, duan2019Grobner}.

Paulino et al. investigated a Bayesian binary regression model with misclassified outcomes \citep{paulino2003binomial}. The problem they considered is a binary regression model of human papillomavirus (HPV) status with three binary risk factors, where the data have $2^3-1=7$ degrees of freedom, yet the model has $2+4=6$ parameters to be estimated. Although the data degrees of freedom are larger than the number of parameters to be estimated, it was observed by the authors that ``model nonidentifiability requires special care in the prior elicitation process because, as is well known, posterior inference for  nonidentifiable parameters is strongly influenced by the prior, even for increasingly large sample sizes." In the frequentist literature, a similar case was studied by \cite{jones2010identifiability}, in which two correlated diagnosis tests were evaluated for three different populations in the absence of a gold standard test. The data degrees of freedom equal the number of parameters yet the model is not locally identifiable and therefore not globally identifiable.

Even though sufficient conditions for local identifiability in some related settings have been established  \citep{jones2010identifiability}, the conditions for local identifiability do not necessarily imply the model is globally identifiable. In fact, the study of global identifiability boils down to solving a system of multivariate nonlinear equations, which is often a challenging mathematical problem and has to be tackled case by case. See \cite{drton2011global} for an example on global identifiability of acyclic mixed graphical models. Recently, \cite{duan2019Grobner} verified that several examples in \cite{jones2010identifiability} are indeed globally identifiably by applying the Gr\"obner basis approach  \citep{sturmfels2005grobner,drton2008lectures}, a powerful computational algebraic tool for solving polynomial equations. However, such an approach is not generally applicable because the nonlinear equations in our setting may not be written as polynomial equations even after reparametrization. Even if it is feasible, the conditions are often difficult to interpret and may not be necessary for global identifiability. Despite the practical importance of the logistic regression models with surrogate outcomes, the global identifiability of such model is still a largely open problem.

In this paper, we establish the necessary and sufficient conditions for the global identifiability of logistic regression models with misclassification, and study a hypothesis testing procedure of regression coefficients in the misclassified logistic regression model when the model is not identifiable under the null. In particular, our results are applicable to models with multiple covariates, which can be either continuous or discrete. Moreover, our identifiability conditions are both necessary and sufficient, and are easy to interpret. The proof of necessity is mainly achieved by proof by contradiction. To show sufficiency, we develop an approach called submodel analysis, which breaks down the problem to simpler identifiability problems for some submodels. One key step in our main proof is adapted from the classical Picard-Lindel\"{o}f existence theorem in the field of ordinary differential equations. In addition to the identifiability conditions, we propose valid hypothesis testing procedures in the case where the model is not identifiable under the null hypothesis.

\section{Background}

We consider logistic regression models with misclassified outcomes as in \cite{grace2016statistical}. Let $Y$ be a binary outcome and $x$ be the corresponding $p$-dimensional covariates. The logistic regression model assumes
\begin{equation}\label{eq:0}
{f}\left(\textrm{Pr}(Y=1 | x)\right) = \beta_0+x^{T}\beta_x,
\end{equation}
where $f(a)=\log\left\{a/(1-a)\right\}$. Instead of observing $Y$, we observe an error-prone surrogate $S$, which is believed to have misclassification errors, sensitivity $\alpha_1$, defined as $\alpha_1 = P (S=1|Y=1)$,  and specificity $\alpha_2$, defined as $\alpha_2 = (S=0|Y=0)$, being independent of the covariates. Denote the collection of parameters by $\theta=(\alpha_1, \alpha_2, \beta_0, \beta_x)$. The model for the surrogate $S$ is
\begin{equation}\label{eq:1}
\textrm{Pr}(S=1 | x) = (1-\alpha_2)+(\alpha_1+\alpha_2-1)h(\beta_0+x^{T}\beta_x),
\end{equation}
where $h(u)=\expit(u) = \exp(u)/\left\{ 1+\exp(u) \right\}$. Model~(\ref{eq:1}) is globally identifiable if and only if $\textrm{Pr}(S=1|x;\theta)=\textrm{Pr}(S=1|x;\theta')$ for all $x$ almost surely implies $\theta=\theta'$.  

It is relatively well-known that model~(\ref{eq:1}) is not identifiable in the unconstrained parameter space $[0,1]^2\times {\R}^{p+1}$, where at least two sets of parameters, $(\alpha_1, \alpha_2, \beta_0,\beta_x)$ and $(1-\alpha_1, 1-\alpha_2, -\beta_0,-\beta_x)$, dictate the same probability distribution in model~(\ref{eq:1}), a phenomenon known as ``label switching" by re-defining the unobserved outcome $Y$ by $1-Y$. Such label switching problem can be eliminated by imposing the constraint $\alpha_1+\alpha_2>1$. However, it remains unknown under what further conditions model~(\ref{eq:1}) is globally identifiable. Doing so is the topic of Section \ref{sec3}.

\section{Global Identifiability of Misclassified Logistic Regression Models}\label{sec3}

In this section, we establish necessary and sufficient conditions for global identifiability of parameters in logistic regression models with response misclassification. Hereafter, we restrict the parameter space of $\theta$ by imposing $\alpha_1+\alpha_2>1$, i.e., the parameter space of $\theta$ is\\ $\Theta=\left\{[0, 1]^2 \times \R^{p+1}: \alpha_1+\alpha_2>1 \right\}$. In addition, we denote the support of covariates $x$ by $\mathcal{D}$, which is defined as the set of values in $\R^p$ with non-zero probability measure. In particular, we allow the covariates to be continuous or discrete.

By the definition of global identifiability, the goal is to find the necessary and sufficient conditions under which the set of $(p+3)$-dimensional nonlinear equation $\textrm{Pr}(S=1|x;\theta)=\textrm{Pr}(S=1|x;\theta')$ holds for all $x\in\mathcal{D}$ implying $\theta=\theta'$.  To solve this problem in the general case with an arbitrary $p\geq 1$, we propose an approach based on a submodel analysis. The main idea is to break down the identifiability problem for model (\ref{eq:1}) with $p$-dimensional covariates to several simpler identifiability problems for some well designed submodels with only $p=1$ or $p=2$ covariates by fixing some of the entries of $x$ at 0. If $0$ is not contained in the support of $x$, some careful reparametrization is crucial. The detailed explanation of the submodel analysis is given in Remark \ref{rem_sub}. To apply the submodel analysis, the following two lemmas on  model identifiability with $p=1$ or $2$ covariates are the fundamental building blocks.

\begin{lemma}\label{thm1}
	For a univariate logistic regression model with unknown misclassification probabilities, i.e., the setting of $p=1$ in model~(\ref{eq:1}), the model parameter $\theta$ is globally identifiable in $\Theta$ if and only if both of the following conditions are satisfied:
	\begin{enumerate}
		\item $\mathcal{D}$ contains at least four unique values;
		\item $\beta_x\neq 0$.
	\end{enumerate}
\end{lemma}

A proof of Lemma 1 is shown in Section \ref{app}. The proof of necessity is achieved by proof by contradiction. The sufficiency is proved by investigating the uniqueness of the solution to a system of nonlinear equations. The nonlinear equations involve not only the unknown parameters $\theta$ but also four unique points in the support of $x$. Since we would like to cover the case that $x$ is discrete, we cannot simply choose $x=0$ or any arbitrary real number as in the continuous case $x\in R$. This further increases the complexity of the equations and the difficulty of investigating the solution set.  One key step in our proof is to lower bound a function of interest in its full domain. To this end, we use the Gronwall inequality which controls the Picard iterative series from one interval to another in the classical Picard-Lindel\"{o}f existence theorem (see Chapter 1 and Chapter 2 of \cite{coddington1955theory}),
we propose to approximate the function on the edge point of a well-constructed
grid and extend the grid little by little iteratively.

Lemma 1 offers interesting insights. For a univariate logistic regression without misclassification, identifiability only requires two unique values of covariate $x$ for the intercept and coefficient parameters to be identifiable.  In contrast, for a univariate logistic regression with misclassification, Lemma 1 states that two additional values in the support of the univariate covariate are needed. In addition, different from logistic regression without misclassification, it further requires the regression coefficient $\beta_x$ to be non-zero to ensure global model identifiability.

Lemma 1 has an important implication in statistical inference. If we are interested in testing the association between the outcome and the covariate, a formal hypothesis testing procedure of $H_0: \beta_x= 0$ would lead to the Davies problem \citep{davies77,davies1987hypothesis}, because all remaining parameters $(\beta_0, \alpha_1, \alpha_2)$ are non-identifiable under the null hypothesis. Interestingly, we show in Theorem \ref{thm4} of Section \ref{sec5} that for testing $\beta_x= 0$, although the remaining parameters are not identifiable, the likelihood ratio test statistic still converges to a simple $\chi_1^2$ distribution. This conclusion can be generalized to any $p\ge 1$.



Next, we consider the case with two covariates $x = (x_1,x_2)$ where the corresponding regression coefficient is $\beta_x = (\beta_{x_1}, \beta_{x_2})$. Here we denote the support of $x_1$ to be $\mathcal{D}_1$ and the support of $x_2$ to be $\mathcal{D}_2$. To avoid the degenerate case, the support of each covariate has to contain at least two values, otherwise the corresponding regression parameter cannot be distinguished from the intercept $\beta_0$.  In addition, we assume the support of $x = (x_1,x_2)$, denoted by $\mathcal{D}$, is the Cartesian product of $\mathcal{D}_1$ and $\mathcal{D}_2$, i.e, $\mathcal{D} = \mathcal{D}_1\times\mathcal{D}_2$.

We have the following identifiability results for model~(\ref{eq:1}) with two covariates.

\begin{lemma}\label{thm2}
	For a logistic regression model of two covariates with unknown misclassification probabilities, i.e., the setting of $p=2$ in model~(\ref{eq:1}), the model parameter $\theta$ is globally identifiable in $\Theta$ if and only if one of the following two conditions is satisfied:
	\begin{enumerate}
		\item there exists $k\in \{1,2\}$, such that $\mathcal{D}_k$ contains at least four unique values, and $\beta_{x_k}\neq 0$;
		\item $\mathcal{D}_1$ or $\mathcal{D}_2$ contains at least three unique values, and both $\beta_{x_1}$ and $\beta_{x_2}$ are non-zero.
	\end{enumerate}
\end{lemma}
A proof of Lemma 2 is provided in the Supplementary Material. From Condition 1, if $\mathcal{D}_1$ contains at least four unique values and $\beta_{x_1}\neq 0$, then the model parameters are identifiable, regardless of whether $\beta_{x_2}$ is zero or not or the number of values in $\mathcal{D}_2$. From Condition 2, if both regression coefficients are non-zero, the assumption on $\mathcal{D}$ can be relaxed, where one of $\mathcal{D}_1$ and $\mathcal{D}_2$ only needs to contain no less than three values. This assumption ensures that the data  degrees of freedom is no less than the number of parameters.

\begin{remark}[Submodel analysis]\label{rem_sub}
	When studying identifiability conditions for models with multiple covariates, we develop a new approach called the submodel analysis, which is elaborated as follows. As an illustration, we focus on condition 1 in Lemma 2. Through a proper reparameterization, we can assume that $\mathcal{D}_2$ contains $0$.
	Therefore, by fixing $x_2 = 0$, we obtain a submodel involving only $x_1$. As a direct consequence of Lemma 1, Condition 1 guarantees that the submodel is identifiable, that is the unknown parameters ($\alpha_1, \alpha_2, \beta_0, \beta_{x_1}$) are uniquely identified. Under the assumption that the misclassification parameters $\alpha_1, \alpha_2$ are identifiable, we can further identify $\beta_{x_2}$ by choosing another submodel involving only $x_2$. We call the above procedure the submodel analysis, which greatly simplifies the analysis of identifiability conditions when we have multiple covariates.
\end{remark}

Now we provide necessary and sufficient conditions for the global identifiability of parameters in the general model~(\ref{eq:1}) with $p$ covariates. Denote $x= (x_1, \dots, x_p)$, the corresponding coefficients $\beta_x = (\beta_{x_1},\dots,\beta_{x_p})$ and the supports as $\mathcal{D}_1,\dots,\mathcal{D}_p$. We assume the support for $x$ as $\mathcal{D}$, which can be written as $\mathcal{D} = \mathcal{D}_1\times\dots\times\mathcal{D}_p$.

\begin{thm}\label{thm1}
	The parameter $\theta$ is globally identifiable in $\Theta$ in model~(\ref{eq:1}) if and only if one of the following three conditions holds:
	\begin{enumerate}
		\item there exists $k\in \{1,\dots,p\}$, such that $\mathcal{D}_k$ contains at least four values, and $\beta_{x_k}\neq 0$.	
		\item there exist $j$ and $k \in\{1,\dots,p\}$, such that $\mathcal{D}_j$ or $\mathcal{D}_k$ contains at least three values,  and both $\beta_{x_j}$ and $\beta_{x_k}$ are non-zero.
		\item there exist $i$, $j$ and $k\in\{1,\dots,p\}$, such that $\beta_{x_i}\neq 0$, $\beta_{x_j}\neq 0$, and $\beta_{x_k}\neq 0$.
	\end{enumerate}	
\end{thm}

A detailed proof of Theorem 1 is provided in Section \ref{app}. To the best of our knowledge, this is the first necessary and sufficient conditions for global identifiability of misclassified logistic regression models which applies to both discrete and continuous covariates.

Unlike the logistic regression model where all the regression coefficients can be $0$ (i.e., $Y$ is not associated with $x$), the regression with misspecified binary outcome requires $Y$ must be associated with some of the covariates $x$ for model identifiability. Comparing to Lemma 2, Theorem 1 offers another possibility (i.e., condition 3) for model identifiability. In fact, the condition 3 implies that the misclassified logistic regression model with three binary risk factors considered in \cite{paulino2003binomial} is globally identifiable if and only if all three coefficients are nonzero.

While there exist three possible identifiability conditions, they can be interpreted in the following coherent way. The $p$-dimensional model is globally identifiable if and only if there exist a subset of the $p$ covariates, denoted by $x_A$, where the submodel containing only $x_A$ is identifiable. For the rest of the variables not in $x_A$, we do not have further assumptions on their support or the corresponding coefficients. The three identifiability conditions and the above interpretation are all obtained by the submodel analysis as detailed in Remark \ref{rem_sub}, which is the key technique used in the proof of Theorem 1.



\section{{\color{black} Hypothesis testing of regression coefficients when models are not identifiable}}\label{sec5}

In practice, statistical inference on a subset of the regression coefficients is often of interest. For example, let $x_B$ denote a subset of $x$ with $d\leq p$ variables, and we want to test $H_0: \beta_{x_B} = 0$. If the reduced model is identifiable according to Theorem 1, the likelihood ratio test can be applied and the asymptotic distribution is $\chi^2_d$. However, when the reduced model is not identifiable and $d<p$, the limiting distribution is no longer chi-square distributed; this known as the Davies problem \citep{davies77,davies1987hypothesis}.

More specifically, we consider statistical inference for regression coefficients in the misclassified logistic regression model~(\ref{eq:1}).The log likelihood function for the parameter $\theta$ is
\begin{eqnarray}
\log L(\theta)&=&\sum_{i=1}^{n} \log \left\{ \textrm{Pr}(S_{i} | x_{i}) \right\}
=\sum_{i=1}^{n} \left\{ S_{i}\log p_i^*+(1-S_{i})\log (1-p_i^*) \right\},\nonumber
\end{eqnarray}
where $p_i^*=\textrm{Pr}(S_i=1 | x_i)$ is defined in equation~(\ref{eq:1}).

Denote $\beta_x=(\eta^{T},\gamma^{T})^{T}$, where $\eta$ corresponds to the $d$ dimensional covariates of interest, and $\gamma$, {the remaining $p-d$ dimensional regression coefficients.} We are interested in testing
$$H_0:\eta=0.$$ 

The asymptotic distributions of the test statistics for $H_0$ depend on whether the reduced model under $H_0$ is identifiable. If the reduced model{ is identifiable, }standard asymptotic arguments \citep{cox1979theoretical} yield that the likelihood ratio and score tests converge weakly to $\chi^2_d$. The details are omitted for limited space. If in the reduced model some parameters are not identifiable under $H_0$, the asymptotic behaviors of the likelihood ratio and score tests are known to be non-regular and the asymptotic distributions may not be $\chi^2$ \citep{davies77,liu03, zhu06}.

Letting $\alpha = (\alpha_1, \alpha_2)$, we now construct the score test and study its asymptotic distribution. Let $U_n(\alpha,\eta,\gamma)$ be the score function for $\eta$, i.e.,
\begin{align*}
&U_n(\alpha,0,\gamma)\\&=\sum_{i=1}^n\left\{s_i\frac{(\alpha_1+\alpha_2-1)f'(t_{i2})z_{i1}^{T}}{1-\alpha_2+(\alpha_1+\alpha_2-1)f(t_{i2})}-(1-s_i)\frac{(\alpha_1+\alpha_2-1)f'(t_{i2})z_{i1}^{T}}{\alpha_2-(\alpha_1+\alpha_2-1)f(t_{i2})}\right\},
\end{align*}
where $t_{i2}=z_{i2}^{T}\gamma$. For a given $\alpha$, we further partition the information matrix $I_{\beta\beta}$ for $\beta=(\eta,\gamma)$ as $I_{\eta\eta}$, $I_{\gamma\eta}$ and $I_{\gamma\gamma}$. The detailed expressions of the partitions are provided in the Supplementary Material. Denote $I_{\eta|\gamma}=I_{\eta\eta}-I_{\eta\gamma}I^{-1}_{\gamma\gamma}I_{\gamma\eta}$ be the partial information matrix of $\eta$ in the presence of $\gamma$.
For each fixed $\alpha$, the score test for $H_0$ is given by
$$
T_2(\alpha)=U_n^{T}(\alpha,0,\widetilde{\gamma}(\alpha))\{I_{\eta|\gamma}(\alpha,0,\widetilde{\gamma}(\alpha))\}^{-1}U_n(\alpha,0,\widetilde{\gamma}(\alpha)),
$$
where $\widetilde{\gamma}(\alpha)$ is the maximum likelihood estimator of $\gamma$ in the reduced model for fixed $\alpha$. Following standard Taylor expansion arguments, $T_2(\alpha)$ converges weakly to $\chi^2_d$ under $H_0$ for each fixed $\alpha$. However, calculation of the test statistic $T_2(\alpha)$ requires the specification of $\alpha$ and the power of $T_2(\alpha)$ relies critically on the specified $\alpha$. This creates a practical difficulty and potential loss of power. Alternatively, we could consider the following uniform score test which does not require the specification of $\alpha$, and defined as
$$
T_2=\sup_{\alpha}\{T_2(\alpha)\}.
$$
In general, the uniform score test $T_2$ is more powerful than $T_2(\alpha)$ when the value of $\alpha$ under the alternative is different from the specified $\alpha$ in the test $T_2(\alpha)$. As shown in Theorem \ref{thm4}, $T_2$ is asymptotically equivalent to the likelihood ratio test, defined as
$$
T_1=-2\left[\sup_{\theta\in H_0}\left\{\log L(\theta)\right\}-\sup_{\theta}\left\{\log L(\theta)\right\}\right],
$$
where $\theta = (\alpha, \eta, \gamma)$.

Interestingly, when $\eta = \beta_x$, which means $\gamma$ is empty, the test statistics $T_1$ and $T_2$ have the asymptotic distribution of $\chi^2_d$, the same asymptotic distribution as if they were calculated from a regular model. The results are provided in Theorem~\ref{thm5} and a proof is outlined in Section \ref{pth4}.
\begin{thm}\label{thm5}
	{Assume that the regularity conditions (R1) $\sim$ (R3) in  Section \ref{inference} hold, and $\eta = \beta_x$,} then $T_1$ and $T_2$ both converge weakly to $\chi^2_d$ under $H_0$.
\end{thm}
{When $\eta \ne \beta_x$,} the asymptotic distributions of $T_1$ and $T_2$ are rather complicated because of the nonidentifiability under $H_0$. In general, we have the following results.
\begin{thm}
	\label{thm4}
	Assume that the regularity conditions (R1) $\sim$ (R3) in Section \ref{inference} hold. {When $\eta \ne \beta_x$} under the null hypothesis $H_0$, $T_1$ and $T_2$  converge weakly to $\sup_{\alpha}\{S^{T}(\alpha)S(\alpha)\}$, where $S(\alpha)$ is a $d$-dimensional Gaussian process with variance $I_d$ and autocorrelation matrix $$\textrm{lim}_{n\rightarrow\infty}E[\{I_{\eta|\gamma}(\alpha,0,\gamma)\}^{-1/2}U_n(\alpha,0,\gamma)U_n^{T}(\alpha',0,\gamma)\{I_{\eta|\gamma}(\alpha',0,\gamma)\}^{-1/2}].$$
\end{thm}
A proof is outlined in Section \ref{app}. In practice, the asymptotic distribution in Theorem \ref{thm4} can be approximated by simulation, for example, using the multiplier bootstrap method \citep{chernozhukov2013}: let $(e_i)_{i=1}^n$ be a sequence of i.i.d. standard normal variables independent of the scores. We suppress parameters to denote one resample $W_n$ as:
$$
W_n=\sum_{i=1}^n\left\{\frac{s_i(\alpha_1+\alpha_2-1)f'(t_{i2})z_{i1}^{T}}{1-\alpha_2+(\alpha_1+\alpha_2-1)f(t_{i2})}-\frac{(1-s_i)(\alpha_1+\alpha_2-1)f'(t_{i2})z_{i1}^{T}}{\alpha_2-(\alpha_1+\alpha_2-1)f(t_{i2})}\right\} e_i
$$
{{then by the main result  of \cite{chernozhukov2013} we can obtain the multiplier bootstrap distribution by repeatedly calculating}}
$$
\max_{\text{range of }\alpha}\{W_n^{T}\{I_{\eta|\gamma}\}^{-1}W_n\}.
$$

{In real data, for example, an outcome variable identified from electronic health records (EHR), may be subject to severe unknown sensitivity and specificity, due to data fragmentation \citep{wei2012fragmentation}. True sensitivity verified by physician review to be above 95\% could be actually as low as 70\% if EHR data is collected from a single source or over a short period. To account for this issue, we allow the range of $\alpha$ to be wide, say a grid of $\alpha_1$ and $\alpha_2$ each varying from $0.7$ to $1.0$ by $0.01$. Then essentially we need to deal a high dimension setting with $p = 31^2$ for a practical $n$ of several hundreds. Although the validity of the multiplier bootstrap for maxima of high-dimensional random vector is addressed in \cite{chernozhukov2013}, we suggest practitioners choose the range of $\alpha$ tightly to raise statistical power.}

\section{Discussion}\label{sec6}
Our Theorem \ref{thm1} fully characterizes the global identifiability of the classical logistic regression models with response misclassification. The proof technique borrows ideas from the Picard-Lindel\"{o}f existence theorem, that when combined with the submodel analysis, is, as far as we know, new for studying model global identifiability. With some proper modification, this technique could be potentially applied to the investigation of identifiability conditions of several related models, such as logistic regression with multivariate misclassified outcomes for longitudinal data, and multi-class logistic regression  with misclassified outcomes. To the best of our knowledge, the global identifiability of these related models is still an open problem and warrants further investigation.

The identifiability of both regression coefficients and misclassification parameters relies critically on the parametric assumptions in~(\ref{eq:0}) and~(\ref{eq:1}). The extension of our results to models with nonlinear structures or differential misclassification needs to be studied case by case and is worth future investigation.


\section{Proofs}\label{app}
As we discussed in Section 3, the support of each covariate has to contain at least two values, otherwise the corresponding regression parameter cannot be distinguished from the intercept $\beta_0$. With this condition, we start with an argument that without loss of generality we can assume they all take values ${0,1}$ through shifting and re-scaling each variable.

For example, if $x = (x_1, \dots,x_p)$ takes value in $\mathcal{D} = \mathcal{D}_1\times\mathcal{D}_2\dots\mathcal{D}_p$ where $ \mathcal{D}_j$ contains at least two different values $\{d_{j1}, d_{j2}\}$ for $1\le j\le p$,  assuming in an increasing order. Now we construct a new variable $z = (z_1, \dots,z_p)$, where $z_j =  (x_j-d_{j1})/(d_{j2}-d_{j1})$. In this way each $z_j$ takes the value $0$ or $1$. By reparameterizing the corresponding regression coefficient $\beta_j' = (d_{j2}-d_{j1})\beta_j$ for $1\le j\le p$, and the intercept $ \beta_0' = \beta_0+\sum_{j=1}^{p}d_{j1}\beta_j$, we obtain that the two linear terms satisfy
\[
\beta_0 +\beta^Tx = \beta_0' +\beta'^Tz
\]
for all possible $x \in \mathcal{D}$ and the corresponding $z$.  Since the reparameterization from $(\beta_0, \beta)$ to $(\beta'_0, \beta')$ is a bijection, $\beta_0$ and $\beta$ can be uniquely identified if and only if $\beta'_0$ and $\beta'$ can be uniquely identified.

\subsection{Lemma 1}
In the following, we provide two technical lemmas which assist the proof of Lemma 1. The proofs of the lemmas are provided in the Appendix.

\begin{lem1}
	For some $c_1, c_2, c_3 \in (0, 1)$ and $d>1$, if there exist $(a_0, b_0, \beta_0, \beta_1)$ that satisfies $a_0, b_0\in (0, 1)$, $\beta_1 \ne 1$ and
	\begin{eqnarray*}
		\left\{ \begin{array}{l}
			a_0+b_0\textrm{expit}(\beta_0) = c_1\\
			a_0+b_0\textrm{expit}(\beta_0+\beta_1) = c_2\\
			a_0+b_0\textrm{expit}(\beta_0+d\beta_1) = c_3\\
		\end{array} \right.,
	\end{eqnarray*}
	then for any fixed $a \in (a_0, \min\{c_1,c_2,c_3\})$, there exist a unique solution $\{b(a),\beta_0(a), \beta_1(a)\}$ that satisfies the above system of equation. Moreover, the function $b(a)$ is continuous and differentiable,  $0<b(a)<1$ and $\beta_1(a) \ne 0$ for all $a \in (a_0, \min\{c_1,c_2,c_3\})$.
\end{lem1}

\begin{lem1}
	Define
	\begin{equation*}
	\psi(x,y,d)=x^{d-1}y^{-d}+x^{1-d}y^{d}-(1-d)(x+\f{1}{x})- d(y+\f{1}{y}),
	\end{equation*}
	and
	\begin{equation*}
	\phi(x,y,d)=\psi(x,y,d)\left[   \f{y-x}{x^{d-1}} \Big( \f{y^{d}-x^{d}}{y-x}-d_1x^{d-1}  \Big)          \right]^{-1}
	\end{equation*}
	for some $x, y>0$ and $x\ne y$. We have $\phi(x,y,d)$ is monotonic about $d$ when $d>1$.
\end{lem1}

Now we prove Lemma 1. For clear the illustration, we denote $\beta_1 = \beta_x$.
\textit{Sufficiency:}
We begin by defining the following bijection
\begin{equation*}
a=1-\alpha_1, \quad {\rm and}\quad b=\alpha_1+\alpha_2 -1,
\end{equation*}
and we have $a\in (0,1)$ and $b>0$. Without loss of generality, we assume that $\{0,1,d_1,d_2\}\in \mathcal{D}$. For simplicity here we denote $f(\cdot)=\expit(\cdot)$ which is an increasing function with positive outcome. We can obtain the following set of equations.
The equation set could be written as
\begin{equation}
\left\{
\begin{aligned}
& a+b \expit(\beta_0)=c_1,\\
& a+b \expit(\beta_0+\beta_1)=c_2,\\
& a+b \expit(\beta_0+d_1\beta_1)=c_3,\\
& a+b \expit(\beta_0+d_2\beta_1)=c_4,
\end{aligned}
\right.
\label{equationset1}
\end{equation}
We first focus on the first three equations, i.e.,
\begin{equation}
\left\{
\begin{aligned}
& a+b\expit(\beta_0)=c_1,\\
& a+b \expit(\beta_0+\beta_1)=c_2,\\
& a+b \expit(\beta_0+d_1\beta_1)=c_3.
\end{aligned}
\right.
\label{equationset2}
\end{equation}
Let $(a_0,b_0,\beta_0,\beta_1)$ to be the solution with the smallest $a_0$ among all the other solutions to above set of three equations. From Lemma A.1, we know that for arbitrary fixed $a$ such that $a_0<a < \min\{c_1,c_2, c_3\}$, the equation set (\ref{equationset2}) have an unique solution $(b,\beta_0,\beta_1)$, where we may define $b=g(a)$ for some function $g$ such that (\ref{equationset2}) holds for all $a\in (a_0,\min\{c_1,c_2, c_3\})$ and $b=g(a)$, notice here $g$ is continuous differentiable by implicit function theorem. From the first two equations of (\ref{equationset2}), we have
\begin{equation*}
\left\{
\begin{aligned}
& \beta_0=\logit\Big( \frac{c_1-a}{g(a)}   \Big), \\
& \beta_1= \logit \Big( \frac{c_2-a}{g(a)} \Big) - \logit\Big( \frac{c_1-a}{g(a)}   \Big),
\end{aligned}
\right.
\end{equation*}
and the third equation of (\ref{equationset2}) leads to
\begin{equation*}
a+g(a) \expit \left( (1-d_1)\logit\Big( \frac{c_1-a}{g(a)}  \Big) +d_1 \logit \Big( \frac{c_2-a}{g(a)} \Big) \right)-c_3=0,
\end{equation*}
which holds uniformly for all $a\in(0,\min\{c_1,c_2, c_3\})$. Therefore, we have
\begin{align}
0 &=1+g'(a)\expit(h_1(a)) \label{condition}  \\
&\qquad +g(a)\expit(h_1(a))(1-\expit(h_1(a)))\frac{\partial h_1(a)}{\partial a},  \nn \\
&= 1+ e^{-h_1(a)}+g'(a)+\frac{g(a)}{1+e^{h_1(a)}}\frac{\partial h_1(a)}{\partial a},  \nn
\end{align}
with
\begin{align*}
h_1(a)&=(1-d_1)\logit\Big( \frac{c_1-a}{g(a)}  \Big) +d_1 \logit\Big( \frac{c_2-a}{g(a)} \Big) \\
&= \log\left( \Big(  \frac{c_1-a}{g(a)+a-c_1} \Big)^{1-d_1} \Big(    \frac{c_2-a}{a+g(a)-c_2}   \Big)^{d_1}                    \right).
\end{align*}
Notice that
\begin{align*}
\frac{\partial h_1(a)}{\partial a} = (1-d_1)\frac{-g(a)-(c_1-a)g'(a)}{(c_1-a)(a+g(a)-c_1)}+d_1 \frac{-g(a)-(c_2-a)g'(a)}{(c_2-a)(a+g(a)-c_2)},
\end{align*}
and combine with (\ref{condition}) we get a linear equation of $g'(a)$, i.e.,
\begin{align*}
&1+e^{-h_1(a)}-\frac{g(a)}{1+e^{h_1(a)}} \cdot  \frac{(1-d_1)g(a)}{(c_1-a)(a+g(a)-c_1)}  \\
&-\frac{g(a)}{1+e^{h_1(a)}} \cdot \frac{d_1g(a)}{(c_2-a)(a+g(a)-c_2)} \\
=& g'(a)\left[ -1+\frac{g(a)}{1+e^{h_1(a)}}\frac{(1-d_1)}{a+g(a)-c_1}+\frac{g(a)}{1+e^{h_1(a)}}\frac{d_1}{a+g(a)-c_2}  \right].
\end{align*}
We define $t_1=(c_1-a)/g(a)$, $t_2=(c_2-a)/g(a)$, and by multiplying $1+e^{h_1(a)}$ on the both side, the above equation can be written as
\begin{align*}
&2+\Big( \frac{t_1}{1-t_1} \Big)^{d_1-1}\Big(  \frac{t_2}{1-t_2}   \Big)^{-d_1}+\Big( \frac{t_1}{1-t_1} \Big)^{1-d_1}\Big(  \frac{t_2}{1-t_2}   \Big)^{d_1}\\
&-\frac{(1-d_1)}{t_1(1-t_1)}-\frac{d_1}{t_2(1-t_2)}\\
=& g'(a) \left[ \frac{t_1}{1-t_1}-  \Big( \frac{t_1}{1-t_1} \Big)^{1-d_1}\Big(  \frac{t_2}{1-t_2}   \Big)^{d_1}+d_1\Big( \frac{t_2}{1-t_2}-\frac{t_1}{1-t_1}  \Big)        \right].
\end{align*}
We further let $x=t_1/(1-t_1)$, $y=t_2/(1-t_2)$, which implies $t_1=x/(1+x)$, and
\begin{align}
&2+x^{d_1-1}y^{-d_1}+x^{1-d_1}y^{d_1}-\f{ (1-d_1)(1+x)^2}{x}- \f{d_1(1+y)^2}{y} \label{aaboutb} \\
=& g'(a) \frac{(x-y)}{x^{d_1-1}}\left( \frac{x^{d_1}-y^{d_1}}{x-y}-d_1x^{d_1-1}   \right). \nn
\end{align}
On the one hand, by Taylor expansion and the fact that  $K(x)=x^{d}$ is an increasing function of $x$ when $d>1$, we see that
\begin{equation*}
\frac{(x-y)}{x^{d_1-1}}\left( \frac{x^{d_1}-y^{d_1}}{x-y}-d_1x^{d_1-1}   \right)<0,
\end{equation*}
for $x,y>0$ and $x\neq y$ ($x=y$ would imply $t_1=t_2$ and $c_1=c_2$, which imply $\beta_1=0$). On the other hand, denote
\begin{equation*}
\psi(x,y,d_1)=x^{d_1-1}y^{-d_1}+x^{1-d_1}y^{d_1}-(1-d_1)(x+\f{1}{x})- d_1(y+\f{1}{y}),
\end{equation*}
and we have $\psi(y,y,d_1)=0$. Since $d_1>1$, we have
\begin{align*}
\partial\psi(x,y,d_1)/\partial x= (d_1-1) \frac{1}{x^2} \left[ \Big( \f{x}{y} \Big)^{d_1}-1    \right] - (d_1-1) \left[  \Big( \f{y}{x}  \Big)^{d_1} -1   \right],
\end{align*}
which implies $\partial\psi(x,y,d_1)/\partial x<0$ when $x\in[0,y)$  and $ \partial\psi(x,y,d_1)/\partial x >0$ when $x\in(y,\infty)$. Hence that $\psi(x,y,d_1)> 0$ for all $x,y>0$, $x\neq y$. Combine this with (\ref{aaboutb}) we conclude $g'(a)<0$, i.e., $b$ is decreasing with $a$.
\vskip 1em
Now, we want to investigate  whether the function
\begin{equation*}
f(a)=a+g(a)\expit(h_2(a))-c_4,
\end{equation*}
has a unique solution for the equation $f(a)=0$ for $a\in(0,\min\{c_1, c_2, c_3\})$, where
\begin{align*}
h_2(a)&=(1-d_2)\logit\Big( \frac{c_1-a}{g(a)}  \Big) +d_2 \logit\Big( \frac{c_2-a}{g(a)} \Big) \\
&= \log\left( \Big(  \frac{c_1-a}{g(a)+a-c_1} \Big)^{1-d_2} \Big(    \frac{c_2-a}{a+g(a)-c_2}   \Big)^{d_2}                    \right),
\end{align*}
and
\begin{equation*}
\frac{\partial h_2(a)}{\partial a} = (1-d_2)\frac{-g(a)-(c_1-a)g'(a)}{(c_1-a)(a+g(a)-c_1)}+d_2 \frac{-g(a)-(c_2-a)g'(a)}{(c_2-a)(a+g(a)-c_2)}.
\end{equation*}
Since
\begin{align*}
f'(a)&=\expit(h_2(a))\left[  1+ e^{-h_2(a)}+g'(a)+\frac{g(a)}{1+e^{h_2(a)}}\frac{\partial h_2(a)}{\partial a}    \right],\\
&= \expit(h_2(a))(1-\expit(h_2(a)))\\
& \quad \times \Bigg[  x^{d_2-1}y^{-d_2}+x^{1-d_2}y^{d_2}- (1-d_2)(x+\f{1}{x})- d_2(y+\f{1}{y})    \\
& \qquad \;\; +g'(a)\Big[  \f{y-x}{x^{d_2-1}} \Big( \f{y^{d_2}-x^{d_2}}{y-x}-d_2x^{d_2-1}  \Big)     \Big] \Bigg],
\end{align*}
by defining
\begin{equation*}
\phi(x,y,d)=\psi(x,y,d)\left[   \f{y-x}{x^{d-1}} \Big( \f{y^{d}-x^{d}}{y-x}-dx^{d-1}  \Big)          \right]^{-1},
\end{equation*}
we have $\phi(x,y,d_1)=-g'(a)$ and
\begin{align*}
f'(a) &=  \expit(h_2(a))(1-\expit(h_2(a))) \left[   \f{y-x}{x^{d_2-1}} \Big( \f{y^{d_2}-x^{d_2}}{y-x}-d_2x^{d_2-1}  \Big)          \right] \\
& \quad \times \Big(  \phi(x,y,d_2)-\phi(x,y,d_1)       \Big).
\end{align*}

\noindent From Lemma A.2, we know that $\phi(x,y,d)$ is a monotonic function of $d$ when $d>1$. Therefore, $\phi(x,y,d_2)-\phi(x,y,d_1)$ as well as $f'(a)$ are either positive or negative for all $a \le \min\{c_1, c_2, c_3\}$, which means $f(a)=0$ has only one solution. We therefore complete the proof of the sufficiency.

\textit{Necessity:}
If $\beta_1 = 0$,  for any $x$, the probability mass function for $S$ is
\[
a+b\textrm{expit}(\beta_0).
\]
Since $b>0$, so there exist $\epsilon>0$, such that $b-\epsilon>0$. We can construct
\begin{equation*}
\left\{ \begin{array}{l}
a' = a  \\
b' = b-\epsilon\\
\beta_0' = \textrm{logit}\{{{b}\over {b-\epsilon}}\textrm{expit}(\beta_0)\},
\end{array} \right.
\end{equation*}
such that for any $x$
\[
a+b\textrm{expit}(\beta_0) = a'+b'\textrm{expit}(\beta_0').
\]
Hence the model is not globally identifiable.

If $x$ can only take three values, without loss of generality, we assume them to be $\{0,1,d_1\}$. Then for the true parameter value $(a_0, b_0, \beta_0, \beta_1)$, we let

\begin{equation*}
\left\{
\begin{aligned}
& a_0+b_0\expit(\beta_0)=c_1,\\
& a_0+b_0 \expit(\beta_0+\beta_1)=c_2,\\
& a_0+b_0 \expit(\beta_0+d_1\beta_1)=c_3.
\end{aligned}
\right.
\end{equation*}

Then by Lemma $A.1$, for any fixed $a \in (a_0,\min\{c_1,c_2,c_3\})$,  the equation set have a unique solution $(b (a),\beta_0(a),\beta_1(a)$. Thus, the model is not globally identifiable.
\subsection{Lemma 2}

\textit{Sufficiency of Condition 1}:
We begin with the sufficiency of Condition 1. Without loss of generality, we assume that  $\{0,1,d_1,d_2\} \subseteq \mathcal{D}_1$, $\beta_1\ne 0$ and $\{0,1\}\subseteq \mathcal{D}_2$.  We can construct the following system of equations:
\begin{equation*}
\left\{ \begin{array}{l}
a+b\textrm{expit}(\beta_0) = c_0\\
a+b\textrm{expit}(\beta_0+\beta_1) = c_1\\
a+b\textrm{expit}(\beta_0+d_1\beta_1) = c_2\\
a+b\textrm{expit}(\beta_0+d_2\beta_1) = c_3\\
a+b\textrm{expit}(\beta_0+\beta_2) = c_4
\end{array} \right.
\end{equation*}

From Lemma 1, we know that $(a,b,\beta_0,\beta_1)$ can be solved uniquely from the first four equations. From the last equation we have
\[
\beta_2 = \textrm{logit}((c_4-a)/b)-\beta_0,
\]
which is identified uniquely .\\

\textit{Sufficiency of Condition 2}:
Now we prove the sufficiency of Condition 2, assume $\{x_{11},x_{12},x_{13}\} \subseteq \mathcal{D}_1$, $\{x_{21},x_{22}\}\subseteq \mathcal{D}_2$, where $x_{11}< x_{12}<x_{13}$ and $x_{21}<x_{22}$. The equation set is given by
\begin{numcases}{}	
	a+b \expit(\beta_0+x_{11}\beta_1+x_{21}\beta_2)=c_1, \nn \\
	a+b \expit(\beta_0+x_{12}\beta_1+x_{21}\beta_2)=c_2, \nn \\
	a+b \expit(\beta_0+x_{13}\beta_1+x_{21}\beta_2)=c_3, \nn \\
	a+b \expit(\beta_0+x_{11}\beta_1+x_{22}\beta_2)=c_4,  \nn \\
	a+b \expit(\beta_0+x_{12}\beta_1+x_{22}\beta_2)=c_5,\nn   \\
	a+b \expit(\beta_0+x_{13}\beta_1+x_{22}\beta_2)=c_6, \nn
\end{numcases}
If $c_1<c_4$, i.e., $\beta_2>0$, since $\expit(\cdot)$ is a monotone function, so we have either $c_1<c_2$, which implies $\beta_1>0$ or $c_1>c_2$, which corresponds to $\beta_1<0$. So we define
\begin{equation*}
\left\{
\begin{aligned}
& \widetilde{\beta}_0=\Big[ \beta_0+x_{11}\beta_1+x_{21}\beta_2 \Big]\cdot I( c_1<c_2) + \Big[ \beta_0+x_{13}\beta_1+x_{21}\beta_2 \Big]\cdot I(c_1>c_2), \\
& \widetilde{\beta}_1= \Big[ (x_{12}-x_{11})\beta_1 \Big]\cdot I( c_1<c_2) + \Big[ (x_{12}-x_{13})\beta_1\Big]\cdot I(c_1>c_2) , \\
&  \widetilde{\beta}_2= (x_{22}-x_{21})\beta_2, \\
& d= \frac{x_{13}-x_{11}}{x_{12}-x_{11}},
\end{aligned}
\right.
\end{equation*}
and if $c_1>c_4$, we just replace $(x_{21},x_{22},c_1,c_2)$ by $(x_{22},x_{21},c_4,c_5)$ correspondingly in the above definition. As a overall result, we would have $d>1$, $\widetilde{\beta}_1, \widetilde{\beta}_2> 0$, i.e, if we abuse the use of notation $(\beta,c)$ a little bit, we may see that, with some reparameterization, the above equation set is equivalent to the following equation set
\begin{numcases}{}	
	a+b \expit(\beta_0)=c_1, \label{1}\\
	a+b \expit(\beta_0+\beta_1)=c_2, \label{2}\\
	a+b \expit(\beta_0+d\beta_1)=c_3, \label{3}\\
	a+b \expit(\beta_0+\beta_2)=c_4,  \label{4}\\
	a+b \expit(\beta_0+\beta_1+\beta_2)=c_5,  \label{5}\\
	a+b \expit(\beta_0+d\beta_1+\beta_2)=c_6, \label{6}
\end{numcases}
with $d>1$, $\beta_1, \beta_2> 0$, hence that $c_1<c_2<c_3$, $c_4<c_5<c_6$.

\textit{Case I:  $d\ge 2$.} We first assume that $d\ge 2$. Keep in mind that $a\in[0,1]\cap \min\{c_i,1\leq i\leq 6\}$ and $b>0$. Let $(a_0,b_0,\beta_0,\beta_1,\beta_2)$ to be the solution with the smallest $a_0$ among all the other solutions to the set of equations (\ref{1})$\sim$(\ref{6}). According to Lemma A.1,  from the first three of the above equation, i.e., $(\ref{1})\sim( \ref{3})$, we know that for arbitrary fixed $a$ satisfying $a_0\leq a < \min\{ c_i,1\leq i\leq 3\}$, there is a unique solution $(b,\beta_0,\beta_1)$, where $b=g(a)$ for some continuous differentiable function $g$.

Then from (\ref{1}), (\ref{2}), (\ref{4}) we obtain
\begin{equation*}
\left\{
\begin{aligned}
& \beta_0=\logit\Big( \frac{c_1-a}{g(a)}   \Big), \\
& \beta_1= \logit \Big( \frac{c_2-a}{g(a)} \Big) - \logit\Big( \frac{c_1-a}{g(a)}   \Big), \\
& \beta_2= \logit \Big( \frac{c_4-a}{g(a)} \Big) - \logit\Big( \frac{c_1-a}{g(a)}   \Big), \\
\end{aligned}
\right.
\end{equation*}
and by plugging it into $(\ref{6})$ we obtain the following function
\begin{align*}
f_1(a) &= a-c_6 +g(a) \\
&\qquad \times \expit \left( -d\logit\Big( \frac{c_1-a}{g(a)}  \Big) +d \logit \Big( \frac{c_2-a}{g(a)} \Big)  +\logit \Big( \frac{c_4-a}{g(a)} \Big) \right),
\end{align*}
In order to show the identifiability, we want to investigate whether $f_1(a)$ has a unique solution to $f_1(a)=0$ for $a\in[a_0,\min\{c_i,1\})$ or not. Denote
\begin{align*}
h_1(a)&=-d\logit\Big( \frac{c_1-a}{g(a)}  \Big) + d\logit \Big( \frac{c_2-a}{g(a)} \Big)  +\logit \Big( \frac{c_4-a}{g(a)} \Big) \\
&= \log\left( \Big(  \frac{g(a)+a-c_1}{c_1-a} \Big)^{-d} \Big(    \frac{c_2-a}{a+g(a)-c_2}   \Big)^d   \Big(    \frac{c_4-a}{a+g(a)-c_4}   \Big)                \right),
\end{align*}
thus,
\begin{align*}
\frac{\partial h_1(a)}{\partial a} &= d\frac{g(a)+(c_1-a)g'(a)}{(c_1-a)(a+g(a)-c_1)}- d\frac{g(a)+(c_2-a)g'(a)}{(c_2-a)(a+g(a)-c_2)}\\
&\qquad - \frac{g(a)+(c_4-a)g'(a)}{(c_4-a)(a+g(a)-c_4)}.
\end{align*}
Since
\begin{align*}
f_1'(a)&=\expit(h_1(a))\left[  1+ e^{-h_1(a)}+g'(a)+\frac{g(a)}{1+e^{h(a)}}\frac{\partial h_1(a)}{\partial a}    \right],\\
&= \expit(h_1(a))(1-\expit(h_1(a))) \\
& \quad \times \Bigg[ 2+  \Big( \frac{1-t_1}{t_1} \Big)^d \Big(  \frac{t_2}{1-t_2}   \Big)^d  \Big(  \frac{t_3}{1-t_3}   \Big) +\Big( \frac{t_1}{1-t_1} \Big)^d \Big(  \frac{1-t_2}{t_2}   \Big)^d  \Big(  \frac{1-t_3}{t_3}   \Big)\\
& \qquad \;\; +\frac{d}{t_1(1-t_1)}-\frac{d}{t_2(1-t_2)}-\frac{1}{t_3(1-t_3)}  \\
& \qquad \;\; +g'(a) \Big( 1 +\frac{d}{1-t_1}-\frac{d}{1-t_2}-\frac{1}{1-t_3}     \Big)  \\
&\qquad \;\; +g'(a) \Big( \frac{1-t_1}{t_1} \Big)^d \Big(  \frac{t_2}{1-t_2}   \Big)^d \Big(  \frac{t_3}{1-t_3}   \Big)  \Bigg]\\
&= \expit(h_1(a))(1-\expit(h_1(a)))\\
& \quad \times \Bigg[ \f{x^d }{y^d z} +\f{d}{x}- \f{d}{y} -\f{1}{z}  +(1+g'(a))\Big( \f{y^d z}{x^d}+dx-dy-z    \Big)   \Bigg],
\end{align*}
where $t_1=(c_1-a)/g(a)$, $t_2=(c_2-a)/g(a)$, $t_3=(c_4-a)/g(a)$, and $x=t_1/(1-t_1)$, $y=t_2/(1-t_2)$, $z=t_3/(1-t_3)$. So $t_1=x/(1+x)$ and $x<y$, $x<z$.
Notice that by equation (\ref{aaboutb}) in the proof of Lemma 1, we have
\begin{align}
& x^{d-1}y^{-d}+x^{1-d}y^{d}-(1-d)(x+\f{1}{x})- d(y+\f{1}{y}) \nn \\
=& g'(a) \frac{(x-y)}{x^{d-1}}\left( \frac{x^{d}-y^{d}}{x-y}-dx^{d-1}   \right) \nn \\
\Leftrightarrow \qquad &
\f{x^d-y^d}{x-y}\f{x^{d-2}}{y^d}-d\f{x^{d-2}}{y} =(1+g'(a)) \left( \frac{x^{d}-y^{d}}{x-y}-dx^{d-1}   \right) \label{de5} \\
\Leftrightarrow \qquad &
\f{y^d-x^d}{xy^d}-\f{d}{x}+\f{d}{y} =(1+g'(a)) \left( \frac{y^{d}-x^{d}}{x^{d-1}}+dx-dy  \right) \label{derivative}
\end{align}
Let
\begin{equation*}
\psi(x,y,z,a)=\f{x^d }{y^d z} +\f{d}{x}- \f{d}{y} -\f{1}{z}  +(1+g'(a))\Big( \f{y^d z}{x^d}+dx-dy-z    \Big),
\end{equation*}
and plug in (\ref{derivative}), we have
\begin{align}
\psi(x,y,z,a)&=\Big( \f{x^d }{y^d z} +\f{d}{x}- \f{d}{y} -\f{1}{z} +  \f{y^d-x^d}{xy^d}-\f{d}{x}+\f{d}{y}  \Big) \nn \\
&\qquad +\bigg[  (1+g'(a))\Big( \f{y^d z}{x^d}+dx-dy-z    \Big) \nn \\
&\qquad \qquad  -(1+g'(a)) \Big( \frac{y^{d}-x^{d}}{x^{d-1}}+dx-dy  \Big)\bigg] \nn \\
&= \Big( -\f{1}{z}\f{y^d-x^d}{y^d} + \f{y^d-x^d}{xy^d}  \Big)+(1+g'(a))\Big( z\f{y^d-x^d}{x^d}-\f{y^d-x^d}{x^{d-1}}            \Big) \nn \\
&= \Big( \f{y^d-x^d}{x^{d-1}}  \Big)(\f{z}{x}-1)\left[ 1+g'(a)+\f{x^{d-1}}{zy^d}    \right]. \nn
\end{align}
Therefore, we have
\begin{align}
f_1'(a)= &\expit(h_1(a))(1-\expit(h_1(a))) \psi(x,y,z,a) \label{de} \\
:= & \expit(h_1(a))(1-\expit(h_1(a)))  \Big( \f{y^d-x^d}{x^{d-1}}  \Big)(\f{z}{x}-1)\phi_1(x,y,z,a). \nn
\end{align}
where
\begin{equation*}
\phi_1(x,y,z,a)=1+g'(a)+\f{x^{d-1}}{zy^d}.
\end{equation*}
In the following, we investigate whether $\phi_1$ is always positive or always negative, which would imply the uniqueness of the solution for $f_1(a)=0$.
\vskip 1em

Since $z>x$, we have
\begin{equation*}
\phi_1(x,y,z,a)<1+g'(a)+\f{x^{d-2}}{y^d}:=B_{11},
\end{equation*}
By (\ref{derivative}), we have
\begin{align*}
&\f{y^d-x^d}{xy^d}-\f{d}{x}+\f{d}{y}+\f{y^d-x^d}{xy^d}+d\f{x^{d-1}}{y^d}-d\f{x^{d-2}}{y^{d-1}} \\
=&(1+g'(a)+\f{x^{d-2}}{y^d}) \left( \frac{y^{d}-x^{d}}{x^{d-1}}+dx-dy  \right),
\end{align*}
and therefore,
\begin{align*}
&\f{2(y^d-x^d)}{xy^d}+\Big(\f{x^{d-1}}{y^{d-1}}+1\Big) \Big(\f{d(x-y)}{xy}\Big)=B_{11}\f{(y-x)}{x^{d-1}}\left( \frac{y^{d}-x^{d}}{y-x}-dx^{d-1} \right) \nn \\
\Leftrightarrow \;\; & \f{y^d-x^d}{y-x}-\f{d(y^{d-1}+x^{d-1})}{2}=B_{11}\f{y^d}{2x^{d-2}}\left( \frac{y^{d}-x^{d}}{y-x}-dx^{d-1} \right) \nn.
\end{align*}
Since
\begin{equation*}
\f{y^d}{2x^{d-2}}\left( \frac{y^{d}-x^{d}}{y-x}-dx^{d-1} \right)>0,
\end{equation*}
so the statement $B_{11} \le 0$ is equivalent to
\begin{align}
\f{y^d-x^d}{y-x}\leq \f{d(y^{d-1}+x^{d-1})}{2} {\text{  for any  }}  y>x>0.\nn \end{align}
Define $\zeta(y|x,d):=(y^d-x^d)-\frac{d}{2}(y^{d-1}+x^{d-1})(y-x)$, the problem is converted to showing $\zeta(y|x,d)\leq 0$ for any $y>x>0$, i.e.,
\begin{align}
\zeta(y|x,d) = (1-\f{d}{2})y^d +\f{d}{2} xy^{d-1} - \f{d}{2} x^{d-1}y +(\f{d}{2}-1)x^d\leq 0 , \label{equiv1}
\end{align}
for any $y>x>0$. Since $\zeta(x)=0$, and $\mu_1(t)=t^{d-2}$ is a non-decreasing function of $t$ when $t>0$, so by the mean value theorem, there exists $x\leq \xi\leq y$, such that,
\begin{align}
\zeta(y)'=\frac{\partial \zeta(y)}{\partial y}&= d(1-\f{d}{2})y^{d-1}+\f{d}{2}(d-1)xy^{d-2}-\f{d}{2}x^{d-1} \label{zetade1} \\
&=\f{d(d-1)(y-x)}{2}\left[ \f{y^{d-1}-x^{d-1}}{(d-1)(y-x)}-y^{d-2}\right] \nn \\
&=\f{d(d-1)(y-x)}{2}\left[ \xi^{d-2}-y^{d-2}\right]\leq 0 . \nn
\end{align}
Since $\xi\leq y$, we have $\xi^{d-2}-y^{d-2}\leq 0$. Thus, when $y>x$ we have $\zeta'(y)\leq 0$.  Since $\zeta(x)=0$, we have $\zeta(y)\leq 0$ for $\forall y>x$. Hence (\ref{equiv1}) holds, i.e., $\phi_1<B_{11}\leq 0$, which implies the set of equations  $(\ref{1})\sim( \ref{6})$ has a unique solution when $d\ge 2$.

\textit{Case II: } Now we consider $1<d<2$, and
by define $\widetilde{\beta}_0=\beta_0+\beta_2$, 	 the set of equations (\ref{1})$\sim$(\ref{6}) can be rewrite as
\begin{numcases}{}	
	a+b \expit(\widetilde{\beta}_0-\beta_2)=c_1, \label{01}\\
	a+b \expit(\widetilde{\beta}_0+\beta_1-\beta_2)=c_2, \label{02}\\
	a+b \expit(\widetilde{\beta}_0+d\beta_1-\beta_2)=c_3, \label{03}\\
	a+b \expit(\widetilde{\beta}_0)=c_4,  \label{04}\\
	a+b \expit(\widetilde{\beta}_0+\beta_1)=c_5,  \label{05}\\
	a+b \expit(\widetilde{\beta}_0+d\beta_1)=c_6, \label{06}
\end{numcases}
Similarly as Case I, by Lemma S.4, from the last three of the above equation, i.e., $(\ref{04})\sim( \ref{06})$, we know that for arbitrary fixed $a$ satisfying $a_0\leq a \leq \min\{ c_i,4\le i\le 6\}$, there exists a unique solution $(\widetilde{g}(a),\widetilde{\beta}_0,\beta_1)$ where $\widetilde{g}$ is continuously differentiable. Notice that, $\widetilde{g}$ may not equal $g$ essentially, but at point $a_0$, we have $\widetilde{g}(a_0)\equiv g(a_0)$. Through similar derivation as Lemma 1, we have $\widetilde{g}$ satisfying,
\begin{equation}
\f{\widetilde{y}^d-\widetilde{x}^d}{\widetilde{y}-\widetilde{x}}\f{\widetilde{x}^{d-2}}{\widetilde{y}^d}-d\f{\widetilde{x}^{d-2}}{\widetilde{y}} =(1+\widetilde{g}'(a)) \left( \frac{\widetilde{y}^{d}-\widetilde{x}^{d}}{\widetilde{y}-\widetilde{x}}-d\widetilde{x}^{d-1}   \right), \label{gde}
\end{equation}
where
\begin{equation*}
\widetilde{x}=\frac{\widetilde{t}_1}{1-\widetilde{t}_1}, \;\; \widetilde{y}=\frac{\widetilde{t}_2}{1-\widetilde{t}_2},\;\; {\rm and}\;\;
\widetilde{t}_1=\frac{c_4-a}{\widetilde{g}(a)},\;\;
\widetilde{t}_2=\frac{c_5-a}{\widetilde{g}(a)}.
\end{equation*}
Now, further denote
\begin{equation*}
\widetilde{z}=\frac{\widetilde{t}_3}{1-\widetilde{t}_3}, \;\; {\rm and}\;\;
\widetilde{t}_3=\frac{c_1-a}{\widetilde{g}(a)},
\end{equation*}
we have $\widetilde{t}_1<\widetilde{t}_4<\widetilde{t}_5$, and consequently, $\widetilde{z}<\widetilde{x}<\widetilde{y}$. By solving equations (\ref{01}), (\ref{04}) and (\ref{05}), we have
\begin{equation*}
\left\{
\begin{aligned}
& \widetilde{\beta}_0=\logit\Big( \frac{c_4-a}{\widetilde{g}(a)}   \Big)=\log \widetilde{x}, \\
& \beta_1= \logit \Big( \frac{c_5-a}{\widetilde{g}(a)} \Big) - \logit\Big( \frac{c_4-a}{\widetilde{g}(a)}   \Big)=\log \widetilde{y}-\log \widetilde{x}, \\
& \beta_2= \logit \Big( \frac{c_4-a}{\widetilde{g}(a)} \Big) - \logit\Big( \frac{c_1-a}{\widetilde{g}(a)}   \Big)=\log \widetilde{x}-\log \widetilde{z}, \\
\end{aligned}
\right.
\end{equation*}
and by plugging it into $(\ref{03})$ we get the following function
\begin{equation*}
f_3(a)= a+\widetilde{g}(a) \expit \left( -d\logit\Big( \frac{c_4-a}{\widetilde{g}(a)}  \Big) +d \logit \Big( \frac{c_5-a}{\widetilde{g}(a)} \Big)  +\logit \Big( \frac{c_1-a}{\widetilde{g}(a)} \Big) \right)-c_3.
\end{equation*}
We hereby investigate whether $f_3(a)$ has a unique solution to $f_3(a)=0$ for $a \in [a_0,\min\{c_i,1\le i\le 3\})$, and $1<d<2$, to investigate the identifiability of the equation set (\ref{01})$\sim$(\ref{06}). Denote
\begin{align*}
h_3(a)&=-d\logit\Big( \frac{c_4-a}{\widetilde{g}(a)}  \Big) + d\logit \Big( \frac{c_5-a}{\widetilde{g}(a)} \Big)  +\logit \Big( \frac{c_1-a}{\widetilde{g}(a)} \Big),
\end{align*}
so we have
\begin{align*}
f_3'(a)= & \expit(h_3(a))(1-\expit(h_3(a)))  \Big( \f{\widetilde{y}^d-\widetilde{x}^d}{\widetilde{x}^{d-1}}  \Big)(\f{\widetilde{z}}{\widetilde{x}}-1)\phi_1(\widetilde{x},\widetilde{y},\widetilde{z},a)
\end{align*}
by (\ref{de}). Since
\begin{equation*}
\expit(h_3(a))(1-\expit(h_3(a)))  \Big( \f{\widetilde{y}^d-\widetilde{x}^d}{\widetilde{x}^{d-1}}  \Big)(\f{\widetilde{z}}{\widetilde{x}}-1)<0,
\end{equation*}
in order to proof $f'(a)<0$ for arbitrary $a\in [a_0,\min{\{c_i,1\}}]$, we only need to show
\begin{equation*}
\phi_1(\widetilde{x},\widetilde{y},\widetilde{z},a)= 1+\widetilde{g}'(a)+\f{\widetilde{x}^{d-1}}{\widetilde{z}\widetilde{y}^d} >0.
\end{equation*}

Notice that by the definition of $\widetilde{x},\widetilde{y},\widetilde{z}$ above, we have $\widetilde{z}<\widetilde{x}<\widetilde{y}$, and therefore
\begin{equation*}
\phi_1(\widetilde{x},\widetilde{y},\widetilde{z},a)>1+\widetilde{g}'(a)+\f{\widetilde{x}^{d-2}}{\widetilde{y}^d}:=B_{21}.
\end{equation*}
Also (\ref{gde}) implies that
\begin{equation*}
2\f{\widetilde{y}^d-\widetilde{x}^d}{\widetilde{y}-\widetilde{x}}\f{\widetilde{x}^{d-2}}{\widetilde{y}^d}-d\Big( \f{\widetilde{x}^{d-2}}{\widetilde{y}} + \widetilde{x}^{d-1}\cdot \f{\widetilde{x}^{d-2}}{\widetilde{y}^d}  \Big)  =\Big( 1+\widetilde{g}'(a)+\f{\widetilde{x}^{d-2}}{\widetilde{y}^d}  \Big) \left( \frac{\widetilde{y}^{d}-\widetilde{x}^{d}}{\widetilde{y}-\widetilde{x}}-d\widetilde{x}^{d-1}   \right),
\end{equation*}
which leads to
\begin{align*}
\f{\widetilde{y}^d-\widetilde{x}^d}{\widetilde{y}-\widetilde{x}}-\f{d(\widetilde{y}^{d-1}+\widetilde{x}^{d-1})}{2}=B_{21}\f{\widetilde{y}^d}{2\widetilde{x}^{d-2}} \left( \frac{\widetilde{y}^{d}-\widetilde{x}^{d}}{\widetilde{y}-\widetilde{x}}-d\widetilde{x}^{d-1}   \right) \nn.
\end{align*}
Since
\begin{equation*}
\f{\widetilde{y}^d}{2\widetilde{x}^{d-2}} \left( \frac{\widetilde{y}^{d}-\widetilde{x}^{d}}{\widetilde{y}-\widetilde{x}}-d\widetilde{x}^{d-1}   \right)>0,
\end{equation*}
showing the statement $B_{21}\geq 0$ is equivalent to show
\begin{align}
\f{\widetilde{y}^d-\widetilde{x}^d}{\widetilde{y}-\widetilde{x}} \geq \f{d(\widetilde{y}^{d-1}+\widetilde{x}^{d-1})}{2} \text{ for any } \widetilde{y}>\widetilde{x}>0.\nn
\end{align}
which can be further converted to show
\begin{align}
\zeta(\widetilde{y}|\widetilde{x},d):=(\widetilde{y}^d-\widetilde{x}^d)-\frac{d}{2}(\widetilde{y}^{d-1}+\widetilde{x}^{d-1})(\widetilde{y}-\widetilde{x}) \geq 0 , \nn
\end{align}
for any $\widetilde{y}>\widetilde{x}>0$, or
\begin{align}	
\zeta(\widetilde{y}|\widetilde{x},d) = (1-\f{d}{2})\widetilde{y}^d +\f{d}{2} \widetilde{x}\widetilde{y}^{d-1} - \f{d}{2} \widetilde{x}^{d-1}\widetilde{y} +(\f{d}{2}-1)\widetilde{x}^d         \geq 0, \label{equiv2}
\end{align}
for any $\widetilde{y}>\widetilde{x}>0$. Since $\zeta(\widetilde{x})=0$, and $\mu_1(t)=t^{d-2}$ is decreasing when $t>0$ and $1<d< 2$, we have
\begin{align}
\zeta(\widetilde{y})'=\frac{\partial \zeta(\widetilde{y})}{\partial \widetilde{y}}&= d(1-\f{d}{2})\widetilde{y}^{d-1}+\f{d}{2}(d-1)\widetilde{x}\widetilde{y}^{d-2}-\f{d}{2}\widetilde{x}^{d-1} \label{zetade2} \\
&= \f{d(\widetilde{y}-\widetilde{x})}{2}\left[ [1+(1-d)]\f{\widetilde{y}^{d-1}}{\widetilde{y}-\widetilde{x}}+(d-1)\f{\widetilde{x}\widetilde{y}^{d-2}}{\widetilde{y}-\widetilde{x}}-\f{\widetilde{x}^{d-1}}{\widetilde{y}-\widetilde{x}}            \right] \nn \\
&=\f{d(\widetilde{y}-\widetilde{x})}{2}\left[ \f{\widetilde{y}^{d-1}-\widetilde{x}^{d-1}}{\widetilde{y}-\widetilde{x}}-(d-1)\f{\widetilde{y}^{d-1}-\widetilde{x}\widetilde{y}^{d-2}}{\widetilde{y}-\widetilde{x}} \right]\nn \\
&=\f{d(\widetilde{y}-\widetilde{x})}{2}\left[ \f{\widetilde{y}^{d-1}-\widetilde{x}^{d-1}}{\widetilde{y}-\widetilde{x}}-(d-1)\widetilde{y}^{d-2}\right]\nn \\
&=\f{d(d-1)(\widetilde{y}-\widetilde{x})}{2}\left[ \widetilde{\xi}^{d-2}-\widetilde{y}^{d-2}\right]\geq 0, \nn
\end{align}
where we used the mean value theorem, as there exist $\widetilde{x}\leq \widetilde{\xi}\leq \widetilde{y}$, s.t.,
\begin{equation*}
\f{\widetilde{y}^{d-1}-\widetilde{x}^{d-1}}{\widetilde{y}-\widetilde{x}}=(d-1)\widetilde{\xi}^{d-2}.
\end{equation*}
Therefore $\zeta'(\widetilde{y})\geq 0$ when $\widetilde{y}>\widetilde{x}$. Since $\zeta(\widetilde{x})=0$, we have $\zeta(\widetilde{y})\geq 0$ for all $ \widetilde{y}>\widetilde{x}$. Hence $\phi_1(\widetilde{x},\widetilde{y},\widetilde{z},a)>B_{21}\geq 0$, and $f'_3(a)<0$, i.e., $f_3(a)$ is decresing and $f_3(a)=0$ has only one solution when $1< d<2$.

Combined Case I and Case II, the sufficiency of Condition 2 is proved. \\
\noindent\textit{Necessity:}
If none of the conditions is satisfied, there are two possible scenarios.
\begin{enumerate}
	\item $\beta_1 = \beta_2 = 0$
	\item Only one of $\beta_1$ and $\beta_2$ is non-zero, and the corresponding covariate variable takes only three values.
\end{enumerate}
If scenario 1 is true, we know $\beta = 0$, and for any $x$, the probability mass function for $S$ is 
\[
a+b\textrm{expit}(\beta_0). 
\]
Since $b>0$, there exist $\epsilon>0$, such that $b-\epsilon>0$. We can construct 
\begin{equation*} 
\left\{ \begin{array}{l} 
a' = a  \\
b' = b-\epsilon\\
\beta_0' = \textrm{logit}\{{{b}\over {b-\epsilon}}\textrm{expit}(\beta_0)\}
\end{array} \right.
\end{equation*}
such that for any $x$ 
\[
a+b\textrm{expit}(\beta_0) = a'+b'\textrm{expit}(\beta_0'). 
\]
Then the model is not globally identifiable.

If scenario $2$ is true, suppose $\beta_1 \ne 0$ and $x_1$ takes only $\{0,1,d\}$. We can construct a set of equations

\begin{equation*}
\left\{
\begin{aligned}
& b\expit(\beta_0)=c_1-a,\\
& b \expit(\beta_0+\beta_1)=c_2-a,\\
& b \expit(\beta_0+d_1\beta_1)=c_3-a.
\end{aligned}
\right.
\end{equation*}

Then by Lemma $A$1, we can state that there exists $a_0<\min\{c_1,c_2,c_3\}$ such that for any fixed $a \in (a_0,\min\{c_1,c_2,c_3\})$,  the equation set have an unique solution $(b,\beta_0,\beta_1)$, where $b=g(a)$ is some continuous differentiable function of $a$. So if we choose $a' = a+\epsilon <\min\{c_1,c_2,c_3\}$ for some $\epsilon>0$, we have $b' = g(a')$, and

\begin{equation*}
\left\{
\begin{aligned}
& \beta_0'=\logit\Big( \frac{c_1-a'}{g(a')}   \Big), \\
& \beta_1'= \logit \Big( \frac{c_2-a'}{g(a')} \Big) - \logit\Big( \frac{c_1-a}{g(a)}   \Big),
\end{aligned}
\right.
\end{equation*}
where $(a',b',\beta_0',\beta_1')$ also satisfy the above system of equations. Thus, the model is not globally identifiable.

\subsection {Theorem 3.1}

\textit{Sufficiency of Condition 1}:\\
Without loss of generality, suppose $k=1$ and $x_1$ takes four values $\{0, 1, d_1, d_2\}$. For the rest of the covariate variables $x_2, \dots, x_p$, we assume that they at least take the value $\{0,1\}$, i.e., $\{0,1\} \in \mathcal{D}_j$ for all $j$. For simplicity of notations, we reparameterize $\alpha_1$ and $\alpha_2$ by
\[
a = 1-\alpha_2, \text{ and } b = \alpha_1+\alpha_2-1.
\]
And the above reparameterization is a bijection, which implies $(a, b)$ is uniquely identified if and only if $(\alpha_1, \alpha_2)$ is uniquely identified. We can then construct the following system of equations containing $p+3$ equations
\begin{eqnarray}
\left\{ \begin{array}{l}
a+b\textrm{expit}(\beta_0) = c_0\\
a+b\textrm{expit}(\beta_0+\beta_1) = c_1\\
a+b\textrm{expit}(\beta_0+d_1\beta_1) = c_2\\
a+b\textrm{expit}(\beta_0+d_2\beta_1) = c_3\\
a+b\textrm{expit}(\beta_0+\beta_2) = c_4\\
\cdots\\
a+b\textrm{expit}(\beta_0+\beta_p) = c_{p+2}	
\end{array} \right.
\end{eqnarray}
By taking $x_j = 0$ for all $2\le j\le p$, we observe that the first four equations are from a reduced model with only one covariate $x_1$ which takes four different values. Therefore, according to Lemma 1, we conclude that from the first four equations, and the condition $\beta_1\ne 0$ we can uniquely identify $(a,b,\beta_0,\beta_1)$.

And for $2\le j\le p$, $\beta_j$ can be uniquely solved as
\[
\beta_j = \text{logit} ({c_{i+2}-a \over b})-\beta_0.
\]
Therefore, the solution of equations (5) is unique and the model is global identifiable.

{
	\textit{Sufficiency of Condition 2}:\\
	Suppose $x_1$ takes $\{0,1,d_1\}$, $\beta_1\ne 0 $, $x_2$ takes $\{0,1\}$ and $\beta_2\ne 0$.
	We can construct the following system of equations
	\begin{eqnarray*}
		\left\{ \begin{array}{l}
			a+b\textrm{expit}(\beta_0) = c_0\\
			a+b\textrm{expit}(\beta_0+\beta_1) = c_1\\
			a+b\textrm{expit}(\beta_0+d_1\beta_1) = c_2\\
			a+b\textrm{expit}(\beta_0+\beta_2) = c_3\\
			a+b\textrm{expit}(\beta_0+\beta_1+\beta_2) = c_4\\
			a+b\textrm{expit}(\beta_0+d_1\beta_1+\beta_2) = c_5\\
			a+b\textrm{expit}(\beta_0+\beta_3) = c_6\\
			\dots\\
			a+b\textrm{expit}(\beta_0+\beta_p) = c_{p+3}
		\end{array} \right.
	\end{eqnarray*}
	By taking $x_j = 0$ for all $3\le j\le p$, we observe that the first six equations are from a reduced model with two covariate $x_1$ and $x_2$, where $x_1$ takes at least three values, and both coefficients are non-zero.  According to Lemma 2, we can uniquely identify the parameters $(a,b,\beta_0,\beta_1, \beta_2)$ from the first five equations.  And for $3\le j\le p$, $\beta_j$  solved uniquely as
	\[can be 
	\beta_j = \text{logit} ({c_{i+3}-a \over b})-\beta_0.
	\]
	Therefore, the model is global identifiable.
}

\textit{Sufficiency of Condition 3:}\\
Without loss of generality, we assume that $\beta_j\ne 0$ for $j = 1,2,3$. Then we have the following system of equations
\begin{eqnarray} \label{eq:11}
\left\{ \begin{array}{l}
a+b\textrm{expit}(\beta_0) = c_0\\
a+b\textrm{expit}(\beta_0+\beta_1) = c_1\\
a+b\textrm{expit}(\beta_0+\beta_2) = c_2\\
a+b\textrm{expit}(\beta_0+\beta_3) = c_3\\
a+b\textrm{expit}(\beta_0+\beta_1+\beta_2) = c_4\\
a+b\textrm{expit}(\beta_0+\beta_1+\beta_3) = c_5\\
a+b\textrm{expit}(\beta_0+\beta_4) = c_6\\
\dots\\
a+b\textrm{expit}(\beta_0+\beta_p) = c_{p+2}
\end{array} \right.
\end{eqnarray}
Using the first four equations of (\ref{eq:11}), we obtain that
\begin{equation}\label{ee}
\left\{ \begin{array}{l}
\beta_0 = \textrm{logit}(\frac{c_0-a}{b})\\
\beta_1 = \textrm{logit}(\frac{c_1-a}{b})-\textrm{logit}(\frac{c_0-a}{b})\\
\beta_2 = \textrm{logit}(\frac{c_2-a}{b})-\textrm{logit}(\frac{c_0-a}{b})\\
\beta_3 =\textrm{logit}(\frac{c_3-a}{b})-\textrm{logit}(\frac{c_0-a}{b})
\end{array} \right.
\end{equation}
Plugging into the fifth and sixth equations we have
\[
a+b\textrm{expit}(\textrm{logit}(\frac{c_1-a}{b})+\textrm{logit}(\frac{c_2-a}{b})-\textrm{logit}(\frac{c_0-a}{b})) = c_4.
\]
and
\[
a+b\textrm{expit}(\textrm{logit}(\frac{c_1-a}{b})+\textrm{logit}(\frac{c_3-a}{b})-\textrm{logit}(\frac{c_0-a}{b})) = c_5.
\]
We can expand $\textrm{expit}(x)$ as $\exp(x)/\{1+\exp(x)\}$, and $\textrm{logit}(x)$ as $\log\{x/(1-x)\}$ in the above function and through some simplification we obtain the following polynomial equations
\begin{equation*}
\left\{ \begin{array}{l}
b c_1 c_2-c_0 c_1 c_2-b c_0 c_4+c_0 c_1 c_4+c_0 c_2 c_4-c_1 c_2 c_4+a^2 (c_0-c_1-c_2+c_4)\\+a (2 c_1 c_2-2 c_0 c_4+b (c_0-c_1-c_2+c_4))=0\\
b c_1 c_3-c_0 c_1 c_3-b c_0 c_5+c_0 c_1 c_5+c_0 c_3 c_5-c_1 c_3 c_5+a^2 (c_0-c_1-c_3+c_5)\\+a (2 c_1 c_3-2 c_0 c_5+b (c_0-c_1-c_3+c_5))=0
\end{array} \right.
\end{equation*}
Solving the above equations we can obtain two solutions for $(a,b)$, which are

Solution 1: $a^{(1)}= -B+\sqrt{B^2-4AC}/(2A)$ and $b^{(1)} = -\sqrt{B^2-4AC}/(2A)$.

Solution 2: $a^{(2)}= -B-\sqrt{B^2-4AC}/(2A)$ and $b^{(2)} = \sqrt{B^2-4AC}/(2A)$,
where $A = c_0c_4 - c_0c_5 - c_1c_2 + c_1c_3 + c_2c_5 - c_3c_4$, $B = c_0c_1c_2 - c_0c_1c_3 - c_0c_1c_4 + c_0c_1c_5 - c_0c_2c_4 + c_0c_3c_5 + c_1c_2c_4 - c_1c_3c_5 + c_2c_3c_4 - c_2c_3c_5 - c_2c_4c_5 + c_3c_4c_5$, and $C = -c_0c_1c_2c_5 + c_0c_1c_3c_4 + c_0c_2c_4c_5 - c_0c_3c_4c_5 - c_1c_2c_3c_4 + c_1c_2c_3c_5$.
Since we constrain so that $\alpha_1+\alpha_2> 1$, the parameter $b$ has to be positive. As a consequence, Solution 1 is not valid. Therefore, Solution 2 is the only solution for $a$ and $b$. Then $\beta_0, \beta_1, \beta_2$ and $\beta_3$ are solved by equations (\ref{ee}), and all  $\beta_j$ where $4\le j \le p$ is identified by
\[
\beta_j = \text{logit} ({c_{i+2}-a \over b})-\beta_0.
\]
Thus the model is globally identified. \\
\textit{Necessity}:
If none of the three condition satisfied, there could be only 3 possible scenarios:
\begin{enumerate}
	\item all $\beta_j =0$ for $j \in \{1, \dots p\}$
	\item there is only one $j \in \{1, \dots p\}$ such that  $\beta_j \ne 0$. And $x_j$ takes less than $4$ values.
	\item there are two non-zero coefficients, $\beta_j$ and $\beta_k$,  $j,k \in \{1, \dots p\}$. And $x_j$, $x_k$ all takes less than $3$ values.
\end{enumerate}
If scenario 1 is true, we know $\beta = 0$, and for any $x$, the probability mass function for $S$ is
\[
a+b\textrm{expit}(\beta_0).
\]
Since $b>0$, there exist $\epsilon>0$, such that $b-\epsilon>0$. We can construct
\begin{equation*}
\left\{ \begin{array}{l}
a' = a  \\
b' = b-\epsilon\\
\beta_0' = \textrm{logit}\{{{b}\over {b-\epsilon}}\textrm{expit}(\beta_0)\}
\end{array} \right.
\end{equation*}
such that for any $x$
\[
a+b\textrm{expit}(\beta_0) = a'+b'\textrm{expit}(\beta_0').
\]
Then the model is not globally identifiable.

If scenario 2 is true, suppose $\beta_1\ne 0$,  $x_1 $ takes only $3$ values, and all $\beta_j=0$ for $j\ne 1$. Then from the necessity of Lemma 1 we know that the model is not identifiable. If  scenario 3 is true, suppose $\beta_1$ and $\beta_2$ are nonzero,  and $x_1$ and $x_2$ only take value $\{0,1\}$. According to the necessity of Lemma 2, the model is not identifiable.

\subsection{Theorems \ref{thm4} and \ref{thm5}}\label{inference}

\subsubsection*{Regularity Conditions}
We consider the following regularity conditions for Theorems \ref{thm4} and \ref{thm5}.

(R1). The parameter $\beta=(\eta^{T}, \gamma^{T})^{T}$ lies in a compact set of a Euclidean space, and the true parameter is an interior point of the space. The misclassification parameters satisfy $0<\alpha_1\leq 1$, $0<\alpha_2\leq 1$ and $\alpha_1+\alpha_2\geq 1+\delta$, for some $\delta>0$.

(R2). The function $E\ell(\alpha,\eta,\gamma)$ has a unique maximizer at the true value of $(\alpha,\eta,\gamma)$.

(R3). The information matrix $I(\theta_0)$ exists and is nonsingular.

\subsubsection*{Proof of Theorem \ref{thm4}}\label{pth4}
\begin{proof}
	For simplicity of notation, assume the dimension of $\eta$ is $d=1$. We consider the following one to one reparametrization, i.e., $(\alpha,\eta,\gamma)\rightarrow(\alpha,\eta,\lambda)$, where $\lambda=\alpha_2-(\alpha_1+\alpha_2-1)f(\gamma)$. Under the new parameterization, the log likelihood under $H_0$ is given by
	\begin{displaymath}
	\sum_{i=1}^n\{s_i\log(1-\lambda)+(1-s_i)\log\lambda\},
	\end{displaymath}
	which does not involve $\alpha$. Hence, $\alpha$ is not identifiable under $H_0$. Note that the maximum likelihood estimator of $\lambda$ under $H_0$ is $\widetilde{\lambda}=1-\bar{s}$, where $\bar{s}=n^{-1}\sum s_i$. One can show that the profile score function for $\eta$ under the new parameterization is
	\begin{displaymath}
	U_n(\alpha,0,\widetilde{\lambda})=\frac{n(\alpha_2-1+\bar{s})(\bar{w}-\bar{z}_1\bar{s})}{\bar{s}(1-\bar{s})},
	\end{displaymath}
	where $\bar{z}_1=n^{-1}\sum z_{i1}$, and $\bar{w}=n^{-1}\sum z_{i1}s_i$. The standardized profile score function is given by
	\begin{displaymath}
	\frac{1}{\sqrt{n}}U_n(\alpha,\widetilde{\lambda})=\frac{\sqrt{n}(\bar{w}-\bar{z}_1\bar{s})(\alpha_2-\lambda)}{\lambda(1-\lambda)}+o_p(1).
	\end{displaymath}
	By the central limit theorem, we have
	\begin{displaymath}
	\frac{\sqrt{n}(\bar{w}-\bar{z}_1\bar{s})}{\nu\sqrt{\lambda(1-\lambda)}}=\frac{\sum_{i=1}^n(z_{i1}-\bar{z}_1)s_i}{\nu\sqrt{n\lambda(1-\lambda)}}\rightarrow N(0,1),
	\end{displaymath}
	where $\nu^2=n^{-1}\sum_{i=1}^n(z_{i1}-\bar{z}_1)^2$ is the sample variance of $z_{i1}$. The partial information is given by $I_{\eta|\lambda}=-n\nu^2(\alpha_2-\lambda)^2/\{\lambda(1-\lambda)\}$. Then
	$$
	T_2(\alpha)=U_n^{T}(\alpha,0,\widetilde{\gamma})(I_{\eta|\gamma})^{-1}U_n(\alpha,0,\widetilde{\gamma})=\bigg(\frac{\sum_{i=1}^n(z_{i1}-\bar{z}_1)s_i}{\nu\sqrt{n\lambda(1-\lambda)}}\bigg)^2+o_p(1).
	$$
	Then, we obtain
	\begin{displaymath}
	T_2=\sup_{\alpha}T_2(\alpha)=\bigg(\frac{\sum_{i=1}^n(z_{i1}-\bar{z}_1)s_i}{\nu\sqrt{n\lambda(1-\lambda)}}\bigg)^2+o_p(1)\rightarrow\chi^2_1.
	\end{displaymath}
	As shown in the proof of Theorem 3, $T_1$ is asymptotically equivalent to $T_2$, which completes the proof.
\end{proof}
\subsubsection*{Proof of Theorem \ref{thm5}}
\begin{proof}
	The information matrix $I_{\beta\beta}(\alpha, 0, \eta)$ is partitioned as follows,
	\begin{align*}
	I_{\eta\eta}=-\sum_{i=1}^n\frac{z_{i1}z_{i1}^{T}}{\{p_i(1-p_i)\}^2}\Big[& \{S_i(1-2p_i)+p_i^2\}(\alpha_1+\alpha_2-1)^2f'(t_{i2})^2\\&-(S_i-p_i)p_i(1-p_i)(\alpha_1+\alpha_2-1)f''(t_{i2})\Big],
	\end{align*}

	\begin{align*}
	I_{\eta\gamma}=-\sum_{i=1}^n\frac{z_{i1}z_{i2}^{T}}{\{p_i(1-p_i)\}^2}\Big[&\{S_i(1-2p_i)+p_i^2\}(\alpha_1+\alpha_2-1)^2f'(t_{i2})^2\\&-(S_i-p_i)p_i(1-p_i)(\alpha_1+\alpha_2-1)f''(t_{i2})\Big],
	\end{align*}
	
	\begin{align*}
	I_{\gamma\gamma}=-\sum_{i=1}^n\frac{z_{i2}z_{i2}^{T}}{\{p_i(1-p_i)\}^2}\Big[&\{S_i(1-2p_i)+p_i^2\}(\alpha_1+\alpha_2-1)^2f'(t_{i2})^2\\&-(S_i-p_i)p_i(1-p_i)(\alpha_1+\alpha_2-1)f''(t_{i2})\Big],
	\end{align*}
	where $t_{i2}=z_{i2}^{T}\gamma$, and  $p_i = \alpha_2-(\alpha_1+\alpha_2-1)f(t_{i2})$. Under the assumptions assumptions (B1) and (B2), the consistency of the maximum likelihood estimator follows from the standard M-estimator theory; see Theorem 5.7 in \cite{van2000asymptotic}. Under the assumption (B3), by Taylor expansions and the law of large number, the profile score function is
	\begin{eqnarray*}
		U_n(\alpha,0,\widetilde{\gamma}(\alpha))&=&U_n(\alpha,0,\gamma)-\frac{\partial^2\ell}{\partial\eta\partial\gamma}\left(\frac{\partial^2\ell}{\partial\gamma\partial\gamma}\right)^{-1}\frac{\partial\ell}{\partial\gamma}+o_p(n^{1/2})\\
		&=&\sum_{i=1}^n \left(\frac{\partial\ell_i}{\partial\eta}-I_{\eta\gamma}I_{\gamma\gamma}^{-1}\frac{\partial\ell_i}{\partial\gamma}\right)+o_p(n^{1/2}),
	\end{eqnarray*}
	where
	$$
	\frac{\partial\ell_i}{\partial\eta}=s_i\frac{(\alpha_1+\alpha_2-1)f'(t_{i2})z_{i1}^{T}}{1-\alpha_2+(\alpha_1+\alpha_2-1)f(t_{i2})}-(1-s_i)\frac{(\alpha_1+\alpha_2-1)f'(t_{i2})z_{i1}^{T}}{\alpha_2-(\alpha_1+\alpha_2-1)f(t_{i2})},
	$$
	$$
	\frac{\partial\ell_i}{\partial\gamma}=s_i\frac{(\alpha_1+\alpha_2-1)f'(t_{i2})z_{i2}^{T}}{1-\alpha_2+(\alpha_1+\alpha_2-1)f(t_{i2})}-(1-s_i)\frac{(\alpha_1+\alpha_2-1)f'(t_{i2})z_{i2}^{T}}{\alpha_2-(\alpha_1+\alpha_2-1)f(t_{i2})}.
	$$
	Since $\{\partial\ell_i/\partial\eta\}$, $\{I_{\eta\gamma}\}$, $\{I_{\gamma\gamma}^{-1}\}$, and $\{\partial\ell_i/\partial\gamma\}$ are all Donsker indexed by $\alpha$ (\cite{van1996weak}, section 2.10). Then, the class of functions $\{\frac{\partial\ell_i}{\partial\eta}-I_{\eta\gamma}I_{\gamma\gamma}^{-1}\frac{\partial\ell_i}{\partial\gamma}:\alpha\}$ is also Donsker. Then $n^{-1/2}U_n(\alpha,0,\widetilde{\gamma}(\alpha))$ converges weakly to a mean zero Gaussian process $U(\alpha)$. Denote by $T_2(\alpha)=S_n^{*T}(\alpha)S_n(\alpha)$, where $S_n(\alpha)=(I_{\eta|\gamma})^{-1/2}U_n(\alpha,0,\widetilde{\gamma}(\alpha))$. Hence, $S_n(\alpha)$ converges weakly to $S(\alpha)=(I_{\eta|\gamma})^{-1/2}U(\alpha)$. By the continuous mapping theorem, $T_2=\sup_\alpha\{S_n^{*T}(\alpha)S_n(\alpha)\}$ converges weakly to $\sup_\alpha\{S^{*T}(\alpha)S(\alpha)\}$, where $S(\alpha)$ is a $d$ dimensional Gaussian process described in the theorem.
	
	Let $\widehat{\eta}(\alpha)$ and $\widehat{\gamma}(\alpha)$ be the MLE of $\eta$ and $\gamma$ for fixed $\alpha$. By Taylor expansions, one can show that
	\begin{eqnarray*}
		T_1&=&\sup_{\alpha}2\{\ell(\alpha,\widehat{\eta}(\alpha),\widehat{\gamma}(\alpha))-\ell(\alpha,0,\widetilde{\lambda}(\alpha))\}\\
		&=&\sup_{\alpha}2[\{\ell(\alpha,\widehat{\eta}(\alpha),\widehat{\gamma}(\alpha))-\ell(\alpha,0,\lambda)\}-\{\ell(\alpha,0,\widetilde{\lambda}(\alpha))-\ell(\alpha,0,\lambda)\}]\\
		&=&\sup_{\alpha}\bigg\{\left(\frac{\partial\ell}{\partial\beta}\right)^{T}\left(-I_{\beta\beta}\right)^{-1}\left(\frac{\partial\ell}{\partial\beta}\right)-\left(\frac{\partial\ell}{\partial\gamma}\right)^{T}\left(-I_{\gamma\gamma}\right)^{-1}\left(\frac{\partial\ell}{\partial\gamma}\right)\bigg\}+o_p(1).
	\end{eqnarray*}
	Using the block matrix inversion formula for $I_{\beta\beta}$, we get
	\begin{eqnarray*}
		T_1&=&\sup_{\alpha}\bigg\{\bigg(\frac{\partial\ell}{\partial\eta}-I_{\eta\gamma}I_{\gamma\gamma}^{-1}\frac{\partial\ell}{\partial\gamma}\bigg)^{T}I_{\eta|\gamma}^{-1}\bigg(\frac{\partial\ell}{\partial\eta}-I_{\eta\gamma}I_{\gamma\gamma}^{-1}\frac{\partial\ell}{\partial\gamma}\bigg)\bigg\}+o_p(1)\\
		&=&\sup_{\alpha}T_2(\alpha)+o_p(1)=T_2+o_p(1).
	\end{eqnarray*}
	Therefore, $T_1$ is asymptotically equivalent to $T_2$. Hence, this completes the proof.
\end{proof}

\setstretch{1.24}
\bibliographystyle{chicago}
\bibliography{logistic}

\end{document}